\documentclass[letterpaper, 10pt]{article}
\usepackage{amsthm,amsmath,amsfonts}

\newcommand{\whp}{whp}

\newcommand{\Undecided}{\textsc{Undecided}}
\newcommand{\Unassigned}{\textsc{Unassigned}}
\newcommand{\AK}{\mathsf{Alon-Kahale}}
\newcommand{\Gal}{\mathsf{Gallager}}
\newcommand{\core}[1]{\mathsf{core}\left(#1\right)}

\newcommand{\PlantedDistCol}[1]{{\cal{G}}^{{\rm plant}}_{n,p,#1}}

\newcommand{\argmin}[2]{\textrm{argmin}_{#1}#2}

\newtheorem{thm}{Theorem}[section]

\newtheorem{prop}[thm]{Proposition}
\newtheorem{defn}[thm]{Definition}
\newtheorem{rem}[thm]{Remark}
\numberwithin{equation}{section}

\title{Message passing for the coloring problem:\\Gallager meets Alon and
Kahale}
\author{Sonny Ben-Shimon \and Dan Vilenchik}
\begin{document}
\maketitle
\begin{abstract}
Message passing algorithms are popular in many combinatorial
optimization problems. For example, experimental results show that
{\em survey propagation} (a certain message passing algorithm) is
effective in finding proper $k$-colorings of random graphs in the
near-threshold regime. In 1962 Gallager introduced the concept of
Low Density Parity Check (LDPC) codes, and suggested a simple
decoding algorithm based on message passing. In 1994 Alon and Kahale
exhibited a coloring algorithm and proved its usefulness for finding
a $k$-coloring of graphs drawn from a certain planted-solution
distribution over $k$-colorable graphs. In this work we show an
interpretation of Alon and Kahale's coloring algorithm in light of
Gallager's decoding algorithm, thus showing a connection between the
two problems - coloring and decoding. This also provides a rigorous
evidence for the usefulness of the message passing paradigm for the
graph coloring problem. Our techniques can be applied to several
other combinatorial optimization problems and networking-related
issues.
\end{abstract}

\section{Introduction and Results}
A \emph{$k$-coloring} $f$ of a graph $G=(V,E)$ is a mapping from
the set of vertices $V$ to $\{1,2,...,k\}$. $f$ is a \emph{proper
coloring} of $G$ if for every edge $(u,v)\in E$, $f(u)\neq f(v)$.
In the graph coloring problem we are given a graph $G=(V,E)$ and
are asked to produce a proper $k$-coloring $f$ with a minimal
possible $k$. The minimal value of $k$ is called the
\emph{chromatic number} of the graph $G$, commonly denoted by
$\chi(G)$. The problem of properly $k$-coloring a $k$-colorable
graph is notoriously hard. It is one of the most famous
NP-Complete problems, and even approximating the chromatic number
within an acceptable ratio is NP-hard \cite{FeigeKilian98}.

\medskip

Message passing algorithms are popular in many combinatorial
optimization problems. For example, experimental results show that
{\em survey propagation} (a certain message passing algorithm) is
effective in finding proper $k$-colorings of random graphs
in the near-threshold regime \cite{SurveyPropagationCol}. In his
seminal work~\cite{Gallager62}, Gallager in 1962 introduced the
concept of Low Density Parity Check (LDPC) codes, and suggested a
simple decoding algorithm based on message passing.

In 1994 Alon and Kahale~\cite{AlonKahale97} introduced, in another
innovative work, a spectral algorithm for coloring sparse (``low
density'') $k$-colorable graphs, and showed that their algorithm
works $\whp$\footnote{By writing $\whp$ we mean with high
probability, i.e., with probability tending to $1$ as $n$ goes to
infinity.} over a certain planted-solution distribution over
$k$-colorable graphs on $n$ vertices.

\subsection{Our Contribution}

In this paper we connect the two latter results. Specifically, we
show that Alon and Kahale's coloring algorithm implicitly contains
Gallager's decoding algorithm, thus asserting an interesting
relation between the two problems -- coloring and decoding.
Theorem \ref{thm:GallagerOnGnpkThm} makes this notion formal.
Hopefully, following this insight, other heuristics that are
useful in decoding algorithms can be applied to the coloring
problem; we are not aware of previous works pinpointing a
connection between the two problems. Our result also gives a
rigorous analysis of a message passing algorithm in the context of
the colorability problem, while up until now only experimental
results showed the usefulness of the message passing paradigm.

\medskip

Our result works in the other direction as well. LDPC codes drew the
attention of many researchers ever since introduced by Gallager in
1962. The problem of graph coloring can be viewed as a specific
instance of the decoding problem -- the codeword being a proper
$k$-coloring of the graph (assuming there is a single one).
Therefore, the analysis of Gallager's algorithm that we provide here
can be useful for the analysis of Gallager's algorithm on random
LDPC codes. Specifically, the study of LDPC codes is concerned with
two main problems. The first is designing ``good" LDPC codes based
on expander graphs and studying their decoding efficiency close to
the code's rate (which usually depends on the expansion properties
of the graph), for example~\cite{RSU01,LMS98}. These works show that
whenever the expander is ``good", one can decode a noisy codeword of
length $n$ close to the code's rate, leaving at most $o(n)$ errors.
The second problem concerns the efficiency of decoding a noisy
codeword with error-ratio some fraction of the code's rate, for
example~\cite{SS96,LMSS98}. In these works it is shown that whenever
the expander is ``good" and the codeword doesn't have too many
errors then typically one can reconstruct the entire codeword (for
example, the rate-1/2 regular codes of Gallager can provably correct
up to 5.17\% errors, or the irregular codes in \cite{LMSS98} with
6.27\%).

Our result concerns the second problem. In our setting, the input
graph $G$ (sampled from the planted distribution over $k$-colorable
graphs, to be defined shortly) on which we analyze Gallager's
algorithm is typically \emph{not} an expander, though it contains a
large subgraph which is an expander. Loosely speaking, we are able
to show that typically Gallager's algorithm converges ``quickly" on
$G$ when starting from a noisy coloring which differs from a proper
one on the color of say $n/120$ vertices to a proper $k$-coloring of
$G$ (up to a small number of vertices which remain ``undecided" and
whose coloring can be completed efficiently). Thus we are able to
analyze Gallager's algorithm while relaxing the demand of expansion
for the entire graph, and replacing it with a more modest
requirement. It could be interesting to apply our techniques back to
the analogous LDPC setting, and in particular to random LDPC codes
similar to the ones introduced by Gallager himself.

\subsection{The Planted Graph Coloring Distribution}
The model we analyze was suggested by Ku\v{c}era \cite{Kuc77} as a
model for generating random $k$-colorable graphs with an arbitrary
average vertex-degree, denoted throughout by $\PlantedDistCol{k}$.
First, randomly partition the vertex set $V=\{1,...,n\}$ into $k$
classes $V_{1},...,V_{k}$, of size $n/k$ each. Then, for every $i
\neq j$, include every possible edge connecting a vertex in
$V_{i}$ with a vertex in $V_{j}$ with probability $p=p(n)$.
Throughout, $\varphi^*:V\rightarrow\{1\ldots,k\}$ denotes the
planted coloring, where the corresponding graph will be clear from
context. $\PlantedDistCol{k}$ is also referred to as the planted
distribution, and is the analog of the planted clique, planted
bisection, and planted SAT distributions, studied e.g. in
\cite{AlonKrivSudCliqe,FeigeKraut,TechReport,flaxman}.

\medskip

$\PlantedDistCol{k}$ is similar to LDPC in some sense. Both
constructions are based on random graphs. In codes, the received
corrupted codeword provides noisy information on a single bit or on
the parity of a small number of bits of the original codeword. In
$\PlantedDistCol{k}$ the adjacency matrix of the graph contains
noisy information about the planted coloring -- embedded in
eigenvectors corresponding to the $(k-1)$  smallest eigenvalues.

\subsection{Our Result}

To avoid a cumbersome presentation we state the results for the case
$k=3$, and point out that the result is easily extended to any fixed
$k$. We call a $k$-coloring $\varphi$ \emph{partial} if some
vertices in $\varphi$ are $\Unassigned$ (that is, are left
uncolored). From here on we will refer to Gallager's decoding
algorithm by $\Gal$ and to Alon and Kahale's coloring algorithm by
$\AK$.

\begin{thm}\label{thm:GallagerOnGnpkThm} Fix $\varepsilon\leq\frac{1}{120}$ and let $G$ be a random graph
in $\PlantedDistCol{3}$, $np \geq C_0$, $C_0=C_0(\varepsilon)$ a
sufficiently large constant. Then running $\Gal$ on $G$ starting
from $\varphi$, some $3$-coloring at distance at most $\varepsilon
n$ from $\varphi^*$, satisfies the following with probability
$1-e^{-\Theta(np)}$:
\begin{enumerate}
  \item  $\Gal(G,\varphi)$ converges to a partial coloring $\varphi'$ after at most $O(\log n)$
  iterations.
  \item  Let $V_A=\{v\in V: \varphi'(v)=\varphi^*(v)\}$, then $|V_A|\geq (1-e^{-\Theta(np)})n$.
  Furthermore, all vertices in $V\setminus V_A$ are $\Unassigned$. \label{i:GallagerOnGnpkThm2}
  \item
  $G[V\setminus V_A]$ can be properly colored (extending the proper coloring of $G[V_A]$) in time $O(n)$.
\end{enumerate}
\end{thm}
\begin{rem}
In item \ref{i:GallagerOnGnpkThm2}, if $np\geq C\log n$ for a
sufficiently large constant $C$, then $\whp$ $V_A=V$, namely all
the vertices are properly colored, and item $3$ becomes void.
\end{rem}
The algorithm $\Gal$ and its adaptation to the coloring problem are
given in Section \ref{sec:GallgersAlgorithm} below. $\AK$ is also
described in Section \ref{sec:AKAlgDescSubs} for the sake of
completeness. For simplicity, we identify  by $\Gal$ both the
decoding algorithm and its adaptation to the colorability problem,
while it will be clear from context which setting is regarded.

\medskip

Comparing our results with the coding setting, known results show
that $\Gal$, applied to carefully designed LDPC codes, typically
recovers the entire codeword (if the noisy codeword has error rate
below the code's rate) or all but $o(1)$-fraction of it (when close
to the rate). In our case, $\Gal$ colors correctly (``decodes") all
but a constant fraction of the vertices, or more precisely, at least
$(1-e^{-\Theta(np)})n$ vertices are properly colored. This is
optimal in some sense as $\whp$ there will be $e^{-\Theta(np)}n$
vertices in $G$ with degree at most $k-2$, each such vertex can
actually take more than one color while still respecting the planted
partition. Furthermore, the vertices which are left $\Unassigned$
(the ones which $\Gal$ fails to ``decode") induce a ``simple" graph
(in the sense that it can be exhaustively 3-colored while respecting
the already-colored part). Nevertheless, in both the decoding
setting and the coloring one it is shown that the quality of
decoding (coloring) improves exponentially with the number of
iterations.


\medskip
Our result is also comparable with the work of \cite{WP} where the
Warning Propagation message passing algorithm was analyzed for the
satisfiability problem. Though Warning Propagation cannot be viewed
as an adaptation of Gallager's algorithm to the satisfiability
problem, the two algorithms are similar in the sense that they both
involve discrete and rather simple (natural) messages. The result in
\cite{WP} is very similar in nature to the one we obtain, though the
initial satisfying assignment can be a random one (due to the
inherent break of symmetry that SAT possesses -- variables can
appear positively or negated).

\subsection{Structure of Paper}
The remainder of the paper is structured as follows. In Section
\ref{sec:GallgersAlgorithm} we present Gallager's decoding algorithm
-- $\Gal$ -- and its adaptation to the coloring problem. We also
present Alon and Kahale's coloring algorithm and show how it
contains $\Gal$, thus justifying the title of this paper. In Section
\ref{sec:PropOfRandomInst} we discuss some properties that a typical
instance in $\PlantedDistCol{3}$ possesses. These properties become
handy when proving Theorem \ref{thm:GallagerOnGnpkThm} in Section
\ref{sec:ProofOfMainThm}. Concluding remarks are given in Section
\ref{sec:Discussion}.

\section{Gallager's Decoding Algorithm}\label{sec:GallgersAlgorithm}
Before presenting $\Gal$ we give a short introduction to linear codes and their graphical representation.

\subsection{Graphical Representation of Codes}

A linear code $\cal{C}$ of length $n$ and dimension $r$ is defined
by a matrix $H\in \{0,1\}^{r\times n}$; a vector
$c=(c_1,\ldots,c_n)\in\{0,1\}^n$ is a \emph{codeword} ($c\in
\cal{C}$) if $Hc\equiv0\,\left(\text{mod }2\right)$. The matrix $H$
can be viewed as the incidence matrix of some bipartite graph
$\mathcal{B}[H]$ with $n$ nodes on the left,
$\left\{x_1,\ldots,x_n\right\}$, and $r$ nodes on the right,
$\left\{C_1,\ldots,C_r\right\}$. Every left-side node corresponds to
a different bit of the message (we refer to $x_i$ as a the $i^{th}$
\emph{variable}), and every right-side node represents a
\emph{constraint} (we refer to $C_j$ as the $j^{th}$
\emph{constraint}). We add the edge $(x_i,C_j)$ if and only if
$H_{i,j}=1$. Every constraint $C_j$ involves the variables in
$N(C_j)=\{x_i:(x_i,C_j)\in E\left(\mathcal{B}[H]\right)\}$, i.e. the
neighbor set of $C_j$. The code $\cal{C}$ can now be defined
equivalently in terms of the bipartite graph $\mathcal{B}[H]$ in the
following way: $c\in\{0,1\}^n$ belongs to $\cal{C}$ if when
assigning $x_i=c_i$, then for every constraint $C_j$ the
exclusive-or of the bits in $N_{\mathcal{B}[H]}(C_j)$ is $0$. If
every variable appears in exactly $s$ constraints, and every
constraint involves exactly $t$ variables, then the bipartite graph
is called $(s,t)$-regular (and so is the corresponding code). If $t$
and $s$ are constants, that is -- every constraint involves at most
a constant number of variables, and every variable takes part in a
constant number of constraints, then the code is called low-density
(an equivalent requirement is that the check-matrix of the code is
sparse) .

\subsection{Gallager's Algorithm}

Gallager's decoding algorithm is a message passing algorithm defined
on the bipartite graph $\mathcal{B}[H]$, which is also referred to
as the \emph{factor graph} in the context of message passing
algorithms. Gallager's algorithm is an example of hard decision
decoding, which signifies the fact that at each step the messages
are derived from local decisions of whether each bit is 0 or 1, and
this is all the information the message contains (as opposed to more
detailed probabilistic information). We note that Gallager also
proposed a belief propagation type decoding algorithm, which uses a
more complicated message set.

\medskip

There are two types of messages associated with every edge
$(x_i,C_j)$ of $\mathcal{B}[H]$; a message from a constraint $C_j$
to a variable $x_i$ and vice versa. Intuitively, we think of the
message $x_i \rightarrow C_j$ as some sort of a majority vote over
the messages $C_{j'}\to x_i$ (for $j\neq j'$). Similarly, $C_j$
sends $x_i$ the ``preferred'' value for $x_i$ -- if the
exclusive-or, $\oplus$, without $x_i$ is $b\in\{0,1\}$, then to keep
the constraint satisfied, $x_i$ should take the value $b$.

Formally, let $\tau\geq 0$ be some fixed integer and $\alpha=(\alpha_1,\alpha_2,...,\alpha_n)$ the received
codeword. Set $\mathcal{C}^{i,b}_j=\{C_{j'}: C_{j'} \rightarrow x_i=b\wedge j'\neq j\}$
for $b\in\{0,1\}$.

\begin{figure}[h]
\framebox[.97\textwidth]{
\parbox{.88\textwidth}{
$x_i \rightarrow C_j = \left\{%
\begin{array}{ll}
    b, & \left|\mathcal{C}^{i,b}_j\right|\geq \tau \\
    \alpha_i, & \hbox{otherwise} \\
\end{array}%
\right.$\\\\
$C_j \rightarrow x_i = \bigoplus_{x\in N_{\mathcal{B}[H]}(C_j)\setminus \{x_i\}}x\rightarrow C_j$
}} \caption{$\Gal$'s messages} \label{fig:GalMsgs}
\end{figure}

It is convenient to define a third message which is not sent during the algorithm,
but is used to calculate the final decoding. For every $x_i$ define
$$B_i = \left\{%
\begin{array}{ll}
    b, & \hbox{$|\{C: C \rightarrow x_i=b\}|\geq \tau$} \\
    \alpha_i, & \hbox{otherwise} \\
\end{array}%
\right.$$
We are now ready to present $\Gal$.
\begin{figure}[h]
\framebox[.97\textwidth]{
\parbox{.88\textwidth}{
$\Gal(\alpha,\tau)$:
\begin{enumerate}
  \item for every edge  $(x_i,C_j):$ set $x_i\rightarrow C_j = \alpha_i$
  \item repeat until no message changed:
    \begin{enumerate}
    \item update in parallel all messages $C_j\rightarrow x_i$.
    \item update in parallel all messages $x_i\rightarrow C_j$.
    \end{enumerate}
  \item return for every $i,$ $x_i=B_i$.
\end{enumerate}
}} \caption{$\Gal$'s algorithm} \label{fig:galAlg}
\end{figure}

The algorithm is allowed not to terminate (the messages may keep
changing every iteration). If the algorithm does terminate, then we
say that it \emph{converges}. In practice, it is common to make an
a-priori limit on the number of iterations, and return failure if
the algorithm does not converge within the given limit. It should be
noted that $\Gal$ is often described with $\tau$ being the size of
the constraint minus $1$.

\subsection{$\Gal$ for Colorability}
Given a $k$-colorable graph $G$ one can define its factor graph
$\mathcal{B}[G]$. The constraints (right hand side nodes) are
the edges of $G$, and the left-hand nodes, the variables,
are the vertices. For $x_i\in V(G)$, and $C_j\in E(G)$ we
add the edge $(x_i,C_j)$ if $x_i\in C_j$. The vector
$c\in\{1,\ldots,k\}^{|V(G)|}$ is a proper $k$-coloring of the graph
if the exclusive-or of each constraint (edge) is \emph{non}-zero
(where the exclusive-or of two integers is the bitwise exclusive-or
of their binary representation). The factor graph of the coloring
problem is a special case of the LDPC factor graph in which every
constraint involves exactly two variables.

We suggest the following interpretation of the messages in Figure
\ref{fig:GalMsgs} for the coloring setting. Fix some constant $\tau$
(Theorem \ref{thm:GallagerOnGnpkThm} is proven with $\tau=1$), and
set as before $\mathcal{C}^{i,b}_j=\{C_{j'}: C_{j'} \rightarrow
x_i=b\wedge j'\neq j\}$, $b\in\{1,\ldots,k\}$. We let $\argmin{}$
denote the index of the minimal elements of the given set.

\begin{figure}[h]
\framebox[.97\textwidth]{
\parbox{.88\textwidth}{
$x_i \rightarrow C_j = \left\{%
\begin{array}{ll}
   \Undecided, & \exists l,m\hbox{ s.t. $l\neq m$ and }\left|\mathcal{C}^{i,l}_j\right|,
   \left|\mathcal{C}^{i,m}_j\right|<\tau;\\
    l, & l=\argmin{t}{\left|\mathcal{C}^{i,t}_j\right|},
\end{array}%
\right.$\\\\
$C_j \rightarrow x_i = x_k \rightarrow C_j.$
}} \caption{$\Gal$'s messages for colorability} \label{fig:GalMsgsCol}
\end{figure}

Though at first glance the messages may appear different than the
ones in Figure \ref{fig:GalMsgs}, this is only due to the fact that
the coloring setting is not a binary one. $x_i$ sends $C_j$ the
minority vote over the colors it receives from the other constrains
(edges) it appears in (or $\Undecided$ if there weren't enough votes
to begin with). The constraint $C_j$ sends the variable $x_i$ which
color it must not take. To see the similarity to the coding setting,
observe that a congruent set of messages for Gallager's algorithm
would be for a variable to send the least popular bit amongst the
messages it received, and for a constraint to send the value which
the variable should not take. The constraint to variable message,
$C_j \rightarrow x_i$, is the same in both settings since in the
coloring problem the exclusive-or is taken over one item (every
constraint -- edge in this case -- contains exactly two variables),
and therefore it is the item itself.

$\Gal$ can now be used for the coloring problem, with the input
$\alpha=(\varphi(x_1),\ldots,\varphi(x_n))$ being the (not necessarily proper) $k$-coloring
vector of the vertices, and the $B_i$'s
defined by
\begin{equation}\label{MPmsgsEq3}
B_i = \left\{%
\begin{array}{ll}
   $\Undecided$, & \exists l,m\hbox{ s.t. $l\neq m$ and } |\{C_j: C_j \rightarrow x_i=l\}|,|\{C_j: C_j \rightarrow x_i=m\}|<\tau\\
   l, & l=\argmin{t}{\left|\{C_j: C_j \rightarrow x_i=t\}\right|}
\end{array}%
\right.
\end{equation}

\subsection{Alon and Kahale's coloring algorithm}\label{sec:AKAlgDescSubs}
The seminal work of Alon and Kahale \cite{AlonKahale97} paved the
road towards dealing with large constant-degree planted
distributions.  They present an algorithm that $\whp$ $3$-colors a
random graph in $\PlantedDistCol{3}$, where $np\geq C_0$ and $C_0$
is a sufficiently large constant. Since our result refers to their
algorithm we give a rough outline of it here, and refer the
interested reader to \cite{AlonKahale97} for the complete details.
\begin{figure*}[!htp]
\begin{center}
\fbox{
\begin{minipage}{\textwidth}
\textsf{Alon-Kahale}($G,k$):\\
\emph{step 1: spectral approximation.}\\
1. obtain an initial $k$-coloring of the graph using spectral methods.\\
\emph{step 2: recoloring procedure.}\\
2. for $i=1$ to $\log n$ do:\\
$\text{ }\text{ }$2.a for all $v \in V$ simultaneously color
$v$ with the least popular color amongst its neighbors.\\
\emph{step 3: uncoloring procedure.}\\
3. while $\exists v \in V$ with less than $np/10$
neighbors colored in some other color do:\\
$\text{ }\text{ }$3.a uncolor $v$.\\
\emph{step 4: Exhaustive Search.}\\
4. let $U \subseteq V$ be the set of uncolored vertices.\\
5. consider the graph $G[U]$.\\
$\text{ }\text{ }$5.a if there exists a connected component of size at least $\log n$ - fail.\\
$\text{ }\text{ }$5.b otherwise, exhaustively extend the coloring
of $V\setminus U$ to $G[U]$.
 \end{minipage}
 }\end{center} \caption{Alon and Kahale's coloring algorithm} \label{Fig:AK97}
 \end{figure*}

The algorithm is composed of three main steps. First using spectral
techniques one typically obtains a $k$-coloring that agrees with the
planted one on many vertices (say $0.99n$). Then follows a refining
procedure (step 2,3) which terminates typically with a partial
$k$-coloring that coincides with the planted coloring on the colored
vertices. This coloring typically colors all but $e^{-\Theta(np)}n$
vertices. Finally, in step 4, the small graph induced by the
uncolored vertices is exhaustively searched for a $k$-coloring that
extends the $k$-coloring of the rest of the graph.

\subsection{Gallager meets Alon-Kahale} \label{sec:ProofOfMainThm}
We now show how $\AK$ actually runs $\Gal$. This observation is not
made in~\cite{AlonKahale97}, nor the connection to message passing
whatsoever.

We first claim that one can unify the recoloring and unassignment
steps of $\AK$ (steps 2 and 3) as follows. In the recoloring step,
assign each vertex with the least popular color amongst its
neighbors, or set it $\Undecided$ if it has less than $\tau$
neighbors colored in some color other than its own ($\tau$ is some
fixed integer, say 2). Using the same arguments as in
\cite{AlonKahale97} one can prove that the algorithm performs the
same with the unified recoloring-uncoloring scheme. It is then left
simply to observe that the unified step of the new $\AK$ is exactly
$\Gal$ (with the messages as in Figure \ref{fig:GalMsgsCol}) since
the message $B_i$ (\ref{MPmsgsEq3}) which implies the color of $v_i$
reads exactly this.
\section{Properties of a Random $\PlantedDistCol{3}$
Instance}\label{sec:PropOfRandomInst} In this section we introduce
some properties of a typical graph in $\PlantedDistCol{3}$. These
properties will come in handy when proving Theorem
\ref{thm:GallagerOnGnpkThm} in the next section. Properties of
similar flavor can be found in \cite{AlonKahale97} along with
complete proofs which are omitted here due to space considerations.
Given a subset $U$ of vertices, $e(U)$ denotes the number of edges
spanned by the vertices of $U$.
\medskip

The first property we discuss is discrepancy. Specifically, a random
graph $\whp$ will not contain a small yet unexpectedly dense
subgraph. Formally,

\begin{prop}\label{prop:NoDenseSubgraphsProp}
Let $G$ be a graph distributed according to $\PlantedDistCol{3}$
with $np\geq C_0$, $C_0$ a sufficiently large constant. Then
$\whp$ there exists \emph{no} subgraph of $G$ containing at most
$n/60$ vertices whose average degree is at least $np/15$.
\end{prop}

The next property we discussed was introduced in \cite{AlonKahale97}
and plays a crucial role in the analysis of the algorithm. Loosely
speaking, a vertex $v$ is considered ``safe'' w.r.t. $\varphi^*$ if
$v$ has many neighbors from every color class of $\varphi^*$ other
than its own, and these neighbors are also ``safe'' w.r.t.
$\varphi^*$. $v$ is ``safe'' in the sense that if it converged to a
wrong color (w.r.t. $\varphi^*$) when the algorithm terminates, then
it has many wrongly-colored neighbors w.r.t. $\varphi^*$, and so do
they. This avalanche effect of many wrongly-colored vertices is,
under certain conditions, very unlikely to happen. Therefore,
typically the ``safe'' vertices converge correctly. The following
definition makes this notion formal.

\begin{defn}\label{def:core}
A set of vertices is called a \textbf{core} of $G=(V,E)$, denoted
$\core{G}$ if for every $v\in \core{G}$:
\begin{itemize}
    \item $v$ has at least $np/5$ neighbors in
    $\core{G} \cap V_i$ for every $i\neq \varphi^*(v)$.
    \item $v$ has at most $np/20$ neighbors from $V\setminus
    \core{G}$,
\end{itemize}
\end{defn}

\begin{prop}\label{prop:coreSize}
Let $G$ be a graph distributed according to $\PlantedDistCol{3}$
with $np\geq C_0$, $C_0$ a sufficiently large constant. Then $\whp$
there exists a core $\core{G}$ satisfying
$|\core{G}|\geq(1-e^{-\Theta\left(np\right)})n$.
\end{prop}

The main idea of the proof is to prove that the following procedure
typically outputs a big core. \textsf{Set $X$ to be the set of
vertices having at least $np/4$ neighbors in $G$ in each color class
of $\varphi^*$ other than its own. Then, repeatedly, delete from $X$
any vertex that has less than $np/5$ neighbors in $X$ in some color
class other than its own or more than $np/20$ neighbors not in $X$.}
This procedure clearly defines a core. To see why this core is
typically large -- observe that to begin with very few vertices are
eliminated from the core (since all but
$e^{-\Theta\left(np\right)}n$ vertices have degree, say, $0.99np/3$
in every color class other then their own). If too many vertices
were removed in the iterative step then a small but dense subgraph
exists (as every removed vertex contributes at least $np/20$ edges
to that subgraph). Proposition \ref{prop:NoDenseSubgraphsProp}
bounds the probability of the latter occurring.

\medskip

Next we characterize the structure of the graph induced by the
non-core vertices (the non-core graph).

\begin{prop}\label{prop:SizeOfConnectedCompProp}
Let $G$ be a graph distributed according to $\PlantedDistCol{3}$
with $np\geq C_0$, $C_0$ a sufficiently large constant. Then $\whp$
every connected component in the non-core graph contains $O(\log n)$
vertices.
\end{prop}
\noindent Proposition \ref{prop:SizeOfConnectedCompProp}, whose
complete proof is given in \cite{AlonKahale97}, will not suffice to
prove Theorem \ref{thm:GallagerOnGnpkThm}, and we need a further
characterization of the non-core graph. Using similar techniques to
those used in the proof of Proposition
\ref{prop:SizeOfConnectedCompProp} we prove:
\begin{prop}\label{prop:NoCycle}Let $G$ be a graph distributed according to $\PlantedDistCol{3}$
with $np\geq C_0$, $C_0$ a sufficiently large constant. Then with
probability $1-e^{-\Theta(np)}$ there exists no cycle in the
non-core graph.
\end{prop}
\begin{rem} Observe that Proposition~\ref{prop:NoCycle} is true only
with a \emph{constant} high probability. In fact, $\whp$ the
connected components of the non-core graph are trees with an
additional edge (unicycles). While one may be able to prove
convergence of $\Gal$ on such a structure as well, it will be
considerably more complicated and technical and not too insightful
to the problem that we discuss in this paper. For more details the
reader is referred to \cite{WP}.
\end{rem}

\section{Proof of Theorem \ref{thm:GallagerOnGnpkThm}}\label{s:GallagerOnGnpkPrf}
A graph $G$ is said to be \emph{typical} in $\PlantedDistCol{3}$
if Propositions \ref{prop:NoDenseSubgraphsProp},
\ref{prop:coreSize}, \ref{prop:SizeOfConnectedCompProp} and
\ref{prop:NoCycle} hold for it. The discussion in Section
\ref{sec:PropOfRandomInst} guarantees that $\whp$ $G$ is typical.
Therefore, to prove Theorem \ref{thm:GallagerOnGnpkThm} it
suffices to consider a typical $G$. In the propositions below we
assume that $G$ is such and do not state anew each time.

We will split the proof into two parts corresponding to the three
items of Theorem \ref{thm:GallagerOnGnpkThm}. First we show that
the messages sent between variables and constraints spanned by
$\core{G}$ converge to $\varphi^*$ when starting with a not ``too
noisy'' coloring. Second we show that the messages sent in the
non-core part of the graph converge to $\varphi^*$ or to the
$\Undecided$ message. This combined with Proposition
\ref{prop:NoCycle} implies that the coloring of the $\Undecided$
vertices can be completed in linear time (list coloring of a
tree).

\begin{prop}\label{prop:CoreConverges}
Let $\varphi$ be a $3$-coloring of $G$ that disagree with
$\varphi^*$ on at most $n/120$ vertices. Then $\whp$ all messages of
$\Gal(\varphi,1)$ spanned by the vertices of $\core{G}$ converge
after at most $O\left(\log n\right)$ iterations. Furthermore, the
message $B_i$ of $v_i \in \core{G}$ equals $\varphi^*(v_i)$.
\end{prop}
\begin{proof}
For every vertex $u$ we define $\varphi_{u\to e}^{(i)}$ to be the
value of the message $u\to e$ that the variable (vertex) $u$ sends
the constraint (edge) $e=(u,v)$ in iteration $i$ of $\Gal$
(according to the message defined in Figure~\ref{fig:GalMsgsCol}).
Let $U_i$ be the set of core vertices for which
$\varphi_{(u,e)}^{(i)}\neq \varphi^*(u)$. It suffices to prove that
$|U_i|\leq |U_{i-1}|/2$ (if this is true, then after $\log n$
iterations $U_{\log n}=\emptyset$). Observe that by our assumption
on $\varphi$ -- $|U_0|\leq n/60$. By contradiction, assume that not
in very iteration $|U_i|\leq |U_{i-1}|/2$, and let $j$ be the first
iteration violating the inequality -- $|U_{j}|\geq |U_{j-1}|/2$.
Consider a vertex $u\in U_{j}$. If $\varphi_{u\to e}^{(i)}\neq
\varphi^*(u)$ , then there must be at least $np/20$ vertices $w$
amongst $u$'s neighbors s.t. $\varphi_{w\to e}^{(i-1)}\neq
\varphi^*(w)$. To see this, observe that any vertex in $\core{G}$
has at most $np/20$ neighbors outside $\core{G}$ (by definition of
core), and hence if it has at most $np/20+np/20 = np/10$ wrongly
colored neighbors, then it has at least $np/4-np/10>np/10$ vertices
from the each of the correct colors, and therefore $\varphi_{u\to
e}^{(i)}$ agrees with $\varphi^*(u)$ (contradicting our assumption).
For conclusion, let $U=U_j \cup U_{j-1}$. Then every $u\in U_j$ has
at least $np/10$ neighbors in $U_{j-1}$. Assuming $|U_{j}|\geq
|U_{j-1}|/2$, the average degree in $G[U]$ is at least
$\frac{np/10\cdot|U_{j}|}{|U|}
\geq\frac{np/10\cdot|U_{j}|}{1.5|U_{j}|}\geq np/15$. Recalling that
$|U|\leq 2\cdot n/120\leq n/60$ contradicts Proposition
\ref{prop:NoDenseSubgraphsProp}.
\end{proof}

\begin{prop}\label{p:noncore_convergence}
Let $i_0$ be the iteration in which $\varphi_{u\to
e}^{(i_0)}=\varphi^*(u)$ for every $u\in\core{G}$ ($\varphi_{u\to
e}^{(i)}$ is the message $u \to e$ in the $i^{th}$ iteration). Then
after at most $t=i_0+O(\log n)$ iterations, all messages of
$\Gal(\varphi,1)$ converge. Furthermore, for every $u\in V(G)$,
$\varphi_{u\to e}^{(t)}$ is either $\varphi^*(u)$ or $\Undecided$.
\end{prop}
\begin{proof}
We first note that it suffices to show that all messages leaving
variables will converge as stated, since all messages leaving the
constraints are just repeaters of messages received at the
constraint in the previous iteration. Let
$G_0=G[V\setminus\core{G}]$ and consider the factor graph that
$G_0$ induces (which is also a tree); for a tree edge $(u,e)$ in
the factor graph  we define $\verb"level"(u,e)$ to be $r$ if $r$
is the maximal length of a path between $u$ and a leaf in the
factor graph from which the edge $(u,e)$ is removed. By induction
on the level $r$ we prove that in iteration $r+i_0$ all vertices
at level $r$ transmit either $\Undecided$ or their planted
coloring. Observe that by the assumption in the proposition, at
all iterations $i>i_0$, all messages leaving $\core{G}$ into $G_0$
state the correct color.

The base case is an edge $(u,e)$ at level 0. If
$\verb"level"(u,e)=0$ then $\varphi_{v\to e}^{(i_0)}$ is either
$\varphi^*(v)$ or $\Undecided$ since the only messages that take
part in the calculation are the ones coming from $\core{G}$ (if
any), and they agree with $\varphi^*(u)$. Now consider an edge
$(u,e)$ at level $r$, and consider iteration $i_0+r$. The key
observation is that the level of all edges $(v,e')$, $e'=(u,v)$ is
strictly smaller than that of $(u,e)$ (since when removing $(u,e)$,
the way from $u$ to all the leaves passes through every $(v,e')$).
Now apply the induction hypothesis and repeat the same argument of
the induction step. Finally, since by Proposition
\ref{prop:SizeOfConnectedCompProp} the level of an edge is at most
$O(\log n)$, we have that after $O(\log n)$ iterations (from $i_0$)
$\Gal$ converges to either the proper coloring or to $\Undecided$.
\end{proof}

\noindent This completes the proof of Theorem
\ref{thm:GallagerOnGnpkThm}.

\section{Discussion}\label{sec:Discussion}
Message passing algorithms are a trendy and promising area of
research in many interesting and fundamental optimization problems,
attracting researchers from different disciplines, e.g. statistical
physics, coding theory and more. Some experimental results show the
effectiveness of such algorithms for instances that seem ``hard" for
many other heuristics \cite{SurveyPropagationCol}. Alas, not many of
the message passing algorithms were rigorously analyzed due to the
inherent complexity of this task. Nevertheless, in this work we give
a rigorous analysis of a well-known message passing algorithm in the
context of the colorability problem, pinpointing an interesting
connection between this problem and decoding algorithms, and showing
the possible value of message passing heuristics for colorability.
Our result also works in the other direction, that is, by
introducing new analytical tools to analyze $\Gal$ in the
colorability setting we suggest the same tools for establishing
rigorous results in the ever developing study of LDPC codes.

\medskip

Our results also extend naturally to other partition problems
where spectral techniques provide a first approximation and then
an iterative algorithm is used to find a proper
solution~\cite{AlonKrivSudCliqe,FeigeOfek06,Boppana87}. It would
be interesting to see whether one can adjust $\Gal$ to these
settings and prove similar results to ours.

\medskip

Another interesting question is to study the counterpart of
code-rate in the colorability problem. Specifically, what is the
``farthest" coloring from the planted one, from which one can prove
convergence of $\Gal$? A random 3-coloring will be at distance
roughly $2n/3$ from the planted one. Now if every vertex with $d$
neighbors has exactly $d/3$ neighbors in every color class (which is
what happens in expectation when starting from a random coloring),
then this is a fixed point of the system and $\Gal$ stays ``stuck".
Therefore it is safe to guess $2n/3$ as an upper bound for that
distance. In this work our techniques enable us to prove $n/120$ as
a lower bound, and although no effort was made to optimize the
constants, we do not believe to be able to reach the aforementioned
upper bound. It will be interesting to close this gap.

\subsection*{Acknowledgement} The authors would like to thank Michael
Krivelevich and Simon Litsyn for their helpful comments.

\end{document}